\documentclass[12pt,leqno]{article}
\usepackage{amsmath,amssymb,amscd,latexsym,bm}

\usepackage[all]{xy}

\newcommand{\hoK}{\widehat{\overline{K}}}
\newcommand{\Ker}{\mathrm{Ker}\,}

\newcommand{\Oh}{{\mathcal O}}

\newcommand{\emm}{{\mathfrak{m}}}

\newcommand{\eo}{\mathfrak{o}}

\newtheorem{theorem}{Theorem}%[section]

\newcounter{aufz}
\newcommand{\pr}{\mathrm{pr}}

\makeatletter
\def\moverlay{\mathpalette\mov@rlay}
\def\mov@rlay#1#2{\leavevmode\vtop{%
   \baselineskip\z@skip \lineskiplimit-\maxdimen
   \ialign{\hfil$\m@th#1##$\hfil\cr#2\crcr}}}
\newcommand{\charfusion}[3][\mathord]{
    #1{\ifx#1\mathop\vphantom{#2}\fi
        \mathpalette\mov@rlay{#2\cr#3}
      }
    \ifx#1\mathop\expandafter\displaylimits\fi}
\makeatother

\begin{document}
\title{Invariant Functions on $p$-divisible Groups and the $p$-adic Corona Problem II}
\author{Behzad Nikzad\footnote{Tempered ai, 1920 Yonge St \# 200, Toronto, ON M4S3E2, Canada, \linebreak behzad.nikzad@tempered.ai} \and Christopher Deninger\footnote{Funded by the Deutsche Forschungsgemeinschaft (DFG, German Research Foundation) under Germany's Excellence Strategy EXC 2044--390685587, Mathematics M\"unster: Dynamics--Geometry--Structure and the CRC 1442 Geometry: Deformations and Rigidity, \linebreak c.deninger@uni-muenster.de}}
\date{}
\maketitle
%\centerline{\bf Abstract:} We develop a theory of \'etale parallel transport for vector bundles with numerically flat reduction on a $p$-adic variety. This construction is compatible with natural operations on vector bundles, Galois equivariant and functorial with respect to morphisms of varieties. In particular, it provides a continuous $p$-adic representation of the \'etale fundamental group for every vector bundle with numerically flat reduction. The results in the present paper generalize previous work by the authors on curves.  They can be seen as a $p$-adic analog of higher-dimensional generalizations of the classical Narasimhan-Seshadri correspondence on complex varieties. Moreover, they provide new insights into Faltings' $p$-adic Simpson correspondence between small Higgs bundles and small generalized representations  by establishing a class of vector bundles with vanishing Higgs field giving rise to actual (not only generalized) representations.

%\small 
%~\\[0.3cm]

%\centerline{{\bf 2010 MSC: 14J20, 11G25} } 
%\input{sec1}
\section{Introduction} \label{sec:1}
Consider a $p$-divisible group $G = (G_{\nu})$ over a complete discrete valuation ring $R$ with fraction field $K$ of characteristic zero and perfect residue field $k = R / \pi R$ of positive characteristic $p$. Let $\eo$ be the valuation ring of $C = \hoK$ and set $\eo_n = \eo / \pi^n \eo$ for $n \ge 1$. The group $G_{\nu} (\eo)$ acts by translation on $G_{\nu} \otimes \eo_n$ for all $n$. This note contains a proof of the following result:

\begin{theorem} \label{t1}
Assume that the connected-\'etale sequence for the dual $p$-divisible group $G'$ splits over $\eo$. Then there is an integer $t \ge 1$ such that the cokernel of the natural inclusion into the $G_{\nu} (\eo)$-fixed module 
\[
\eo_n \hookrightarrow \Gamma (G_{\nu} \otimes \eo_n , \Oh)^{G_{\nu} (\eo)}
\]
is annihilated by $p^t$ for all $\nu$  and $n$.
\end{theorem}

Theorem \ref{t1} and variants were shown in \cite{D} under the very restrictive assumption that additionally $\dim G' \le 1$. See \cite{D} \S\,2, \S\,5 for background and motivation regarding the theorem. An unexpected aspect of Theorem \ref{t1} is its relation to the $p$-adic Corona problem. The classical $p$-adic Corona problem was stated by van~der~Put \cite{P}, and solved by him in the one-dimensional case. The solution in the general case is due to Bartenwerfer \cite{B}. 

\textbf{Corona Theorem (van~der~Put, Bartenwerfer)} \textit{Let $C$ be an algebraically closed field with a complete, non-trivial non-archimedean valuation with valuation ring $\eo$. Let $H^{\infty} (\Delta^d) = \eo [[X_1 , \ldots , X_d]] \otimes_{\eo} C$ be the $C$-algebra of bounded analytic functions on the open polydisc $\Delta^d = \{ x \in C^d \mid \| x \| < 1 \}$ where $\| x \| = \max_{1 \le i \le d} |x_i|$. Then for $f_1 , \ldots , f_n \in H^{\infty} (\Delta^d)$ the following conditions are equivalent: \\
(a)  The functions $f_1 , \ldots , f_n$ generate the $C$-algebra $H^{\infty} (\Delta^d)$.\\
(b) There is a constant $\delta > 0$ such that
\[
\max_{1 \le j \le n} |f_j (x)| \ge \delta \quad \text{for all} \; x \in \Delta^d \; .
\]
}

In \cite{D} \S\,4 it was shown that Theorem \ref{t1} follows from the following generalized Corona theorem which is proved in the present note:

\begin{theorem}
\label{t2}
Let $C$ be an algebraically closed field with a complete, non-trivial non-archimedean valuation. For functions $g_1 , \ldots , g_n \in H^{\infty} (\Delta^d)$ the following conditions are equivalent:\\
(a) $(g_1 , \ldots , g_n) \supset (X_1 , \ldots , X_d)$\\
(b) There is a constant $\delta > 0$ such that
\[
\max_{1 \le j \le n} |g_j (x)| \ge \delta \| x \| \quad \text{for all} \; x \in \Delta^d \; .
\]
\end{theorem}

Note that since $(X_1 , \ldots , X_d)$ is maximal in $H^{\infty} (\Delta^d)$, condition (a) is equivalent to $(g_1 , \ldots , g_n) = H^{\infty} (\Delta^d)$ or $(g_1 , \ldots , g_n) = (X_1 , \ldots , X_d)$. For $d = 1$, Theorem \ref{t2} is an immediate consequence of the Corona Theorem as was noted in \cite{D}, Proposition 12. For general $d$, Theorem \ref{t2} was essentially \cite{D} Conjecture 11. Now that it is proved, Theorems 1 and 3 of \cite{D} hold in general and not only for dimensions $d = 1$. In the next section, Theorem \ref{t2} will be deduced from the classical $p$-adic Corona Theorem.

The genesis of this note is a bit unusual. The first author B.N. is an AI expert with a background in mathematics, who recently used ChatGPT to expand ideas c.f. \cite{N} which he developed several years ago while working on his bachelor's thesis. The prompt to look for applications then led AI to a proof of \cite{D} Conjecture 11. The contribution of the second author -- at the suggestion of B.N. -- consisted in checking the proof and writing it up in his own style. While being very pleased and thankful to B.N. that the dimension restrictions in the results of \cite{D} are now removed, the second author is much less impressed by his own failure at the time of writing \cite{D} to see the simple and natural argument that ChatGPT used to deduce Theorem \ref{t2} for all $d \ge 1$ from the $p$-adic Corona Theorem.
\section{Proof of Theorem \ref{t2}} \label{sec:2}

For $g = \sum_{\alpha} c_{\alpha} T^{\alpha}$ in $H^{\infty} (\Delta^d)$ set $\| g\| = \sup_{\alpha} |c_{\alpha}|$. For $x \in \Delta^d$ we have $|g (x)| \le \| g \|$. If (a) holds, we have
\[
X_i = \sum^n_{j=1} a_{ij} g_j \quad \text{for some} \; a_{ij} \in H^{\infty} (\Delta^d) \quad \text{and} \; 1 \le i \le d \; .
\]
Setting $K = \max_{i,j} \| a_{ij} \| > 0$, we get the estimate
\[
\| x \| = \max_i |x_i| \le \max_{i,j} |a_{ij} (x)|  \, |g_j (x)| \le K \max_j |g_j (x)| \quad \text{for} \; x \in \Delta^d \; .
\]
Thus condition (b) follows with $\delta = K^{-1}$. Now assume that (b) holds and that in addition some $g_j$ does not vanish at $x = 0$. Then we have
\[
\max_j |g_j (x)| \ge \delta' > 0 \quad \text{in a neighborhood of} \; 0 \in \Delta^d \; .
\]
Together with (b) we find some $\delta'' > 0$ with 
\[
\max_j |g_j (x)| \ge \delta'' \quad \text{for all} \; x \in \Delta^d \; .
\]
The Corona Theorem now implies that $(g_1 , \ldots , g_n) = H^{\infty} (\Delta^d)$. Thus we may assume that $J = (g_1 , \ldots , g_n) \subset \emm = (X_1 , \ldots , X_d)$ and using (b) we have to show that $J \supset \emm$. This reduction was already made in \cite{D} \S\,4. We can write $g_j = l_j + r_j$ with $l_j$ homogeneous of degree $1$ and $r_j \in \emm^2$. Viewing the $l_j$ as linear maps $l_j : C^d \to C$, we claim that the linear map $L = (l_1 , \ldots , l_n) : C^d \to C^n$ is injective. Otherwise there is a vector $v \in \Ker L$ with $\| v \| = \max |v_i | = 1$. Using that $l_j (tv) = 0$ we find
\[
\max_j |g_j (tv)| = \max_j |r_j (tv)| \le |t|^2 \max_j \| r_j \| \quad \text{for} \; |t| < 1 \; .
\]
On the other hand, condition (b) gives
\[
\max_j |g_j (tv)| \ge \delta |t| \; .
\]
For small $|t|$ these inequalities are incompatible. Hence $L$ is injective and $L^* : (C^n)^* \to (C^d)^*$ therefore surjective. Thus the linear forms $l_1 = L^* (\pr_1) , \ldots , l_n = L^* (\pr_n)$ span $(C^d)^*$ where $\pr_i$ is the $i$-th projection. It follows that the classes of the $g_j$ span $\emm / \emm^2$. Thus $\emm = J + \emm^2$ and therefore $\emm M = M$ for the finitely generated $H^{\infty} (\Delta^d)$-module $\emm / J$. Note that $\emm$ is not contained in the Jacobson radical of $H^{\infty} (\Delta^d)$ since $1 - aX_1$ is not a unit if $|a| > 1$. Nonetheless, the version of the Nakayama Lemma in \cite{M} (1.M) Lemma 1.3 gives an $s \in 1 + \emm$ with $sM = 0$ and hence $s \emm \subset J$. Since $s (0) = 1$, there is some $0 < r < 1$ with $|s (x)| = 1$ for $\| x \| < r$. For $r \le \| x \| < 1$ condition (b) gives $\max_j |g_j (x)| \ge \delta r$. Consequently, 
\[
\max (|s (x)| , |g_1 (x)| , \ldots , |g_n (x)|) \ge \min (1 , \delta r) > 0
\]
for all $x \in \Delta^d$. The Corona Theorem now yields elements $a , b_1 , \ldots , b_n \in H^{\infty} (\Delta^d)$ such that
\[
as + \sum^n_{j=1} b_j g_j = 1 \; .
\]
Using that $s \emm \subset J$ it follows that
\[
X_i = a (s X_i) + \sum^n_{j=1} (b_j X_i) g_j \in J \quad \text{for} \; 1 \le i \le d \; .
\]
Thus $\emm \subset J$ as was to be shown. 
%\input{sec3}
%\input{sec4}
%\bibliographystyle{acm}%alpha
%\bibliography{lit}

\begin{thebibliography}{1}

\bibitem{B}
{\sc Bartenwerfer, W.}
\newblock Die {L}\"osung des nichtarchimedischen {C}orona-{P}roblems f\"ur
  beliebige {D}imension.
\newblock {\em J. Reine Angew. Math. 319\/} (1980), 133--141.

\bibitem{D}
{\sc Deninger, C.}
\newblock Invariant functions on {$p$}-divisible groups and the {$p$}-adic
  corona problem.
\newblock {\em Tokyo J. Math. 33}, 2 (2010), 393--406.

\bibitem{M}
{\sc Matsumura, H.}
\newblock {\em Commutative algebra}, second~ed., vol.~56 of {\em Mathematics
  Lecture Note Series}.
\newblock Benjamin/Cummings Publishing Co., Inc., Reading, MA, 1980.

\bibitem{N}
{\sc Nikzad, B.}
\newblock Anti-isomorphism between submodules of function spaces and separating
  sets.
\newblock https://mathoverflow.net/questions/486809/.

\bibitem{P}
{\sc van~der Put, M.}
\newblock The non-archimedean corona problem.
\newblock In {\em Table {R}onde d'{A}nalyse {N}on {A}rchim\'edienne ({P}aris,
  1972)}, vol.~Tome 102 of {\em Suppl\'ement au Bull. Soc. Math. France}. Soc.
  Math. France, Paris, 1974, pp.~287--317.

\end{thebibliography}

%\input{address}
%\newpage
%\input{sec3}
\end{document}